\documentclass[11pt]{article}

\usepackage{bbm}
\usepackage{euscript}
\usepackage{makeidx}
\usepackage{sec}
\usepackage{pb-diagram}
\usepackage{lamsarrow}
\usepackage{pb-lams}
\usepackage{amsmath}
\usepackage{amsthm}
\usepackage{epsfig}
\usepackage{amssymb}
\usepackage{pifont}
\usepackage{amsbsy}

\theoremstyle{plain}
\newtheorem{theorem}{Theorem}[section]
\newtheorem{lemma}[theorem]{Lemma}

\theoremstyle{definition}

\newtheorem{remark}[theorem]{Remark}

\newcommand\bC{{\mathbb C}}

\newcommand\bZ{{\mathbb Z}}

\DeclareMathOperator{\re}{{\mathfrak R}{\mathfrak e}}

\newcommand{\f}{\EuScript{F}}

\newcommand{\CL}{\EuScript{L}}

\newcommand{\n}{\EuScript{N}}

\newcommand{\p}{\EuScript{P}}
\newcommand{\q}{\EuScript{Q}}

\newcommand{\ra}{\rightarrow}

\def\inpr{\mathbin{\hbox to 6pt{\vrule height0.4pt width5pt depth0pt \kern-.4pt \vrule height6pt width0.4pt depth0pt\hss}}}

\newcommand{\si}{{\sigma}}

\newcommand{\vfi}{{\varphi}}

\newcommand{\pa}{\partial}

\begin{document}

\title{Derangements and asymptotics of Laplace transforms of polynomials.}
\author{Liviu I. Nicolaescu\\Dept. of Mathematics\\
University of Notre Dame\\
Notre Dame, IN 46556-4618\\
nicolaescu.1@nd.edu}

\maketitle

\begin{abstract}
We    describe the behavior as  $n\ra \infty$ of  the Laplace transforms of $P^n$, where  $P$  a fixed complex polynomial. As a consequence  we obtain a new elementary proof  of  an result  of Gillis-Ismail-Offer \cite{GIO} in the combinatorial theory of derangements.
\end{abstract}

\section{Statement of the main results}
\setcounter{equation}{0}

The    generalized derangement problem in combinatorics   can be formulated as follows.  Suppose $X$ is a finite set  and $\sim$ is an equivalence relation on $X$. For each $x\in X$ we denote by $\hat{x}$ the equivalence class of $x$. $\Hat{X}_\sim$ will denote the set of equivalence classes. The counting function of $\sim$ is the function
\[
\nu=\nu_\sim:\hat{X}\ra \bZ,\;\; \nu(\hat{x})=|\hat{x}|.
\]
A $\sim$-\emph{derangement} of $x$ is a permutation $\vfi:X\ra X$ such that
\[
x\not\in \hat{x},\;\;\forall  x\in X.
\]
We denote by $\n(X,\sim)$ the number of  $\sim$-derangements. The ratio
\[
p(X,\sim)=\frac{\n(X,\sim)}{|X|!}
\]
is the probability that a randomly chosen permutation of $X$ is a derangement.

In \cite{EvGil} S. Even and J. Gillis      have  described   a beautiful relationship between  these numbers and  the Laguerre polynomials
\[
L_n(x)={e^x}\frac{d^n}{dx^n}\Bigl(e^{-x} x^{n}\Bigr)=\sum_{k=0}^n\binom{n}{k}\frac{(-x)^k}{k!},\;\;n=0,1,\cdots .
\]
For example
\[
L_0(x)=1,\;\; L_1(x)=1-x,\;\; L_2(x) =\frac{1}{2!}(x^2-4x+2).
\]
We set
\[
L_\sim:=\prod_{c\in \hat{X}}(-1)^{\nu(c)}\nu(c)!L_{\nu(c)}(x).
\]
Observe that the leading coefficient  of $L_\sim$ is $1$. We have the following result.
\begin{theorem}[Even-Gillis]
\begin{equation}
\n(X,\sim)= \int_0^\infty e^{-x} L_\sim(x) dx.
\label{eq: ev-gil}
\end{equation}
\label{th: ev-gil}
\end{theorem}
For a very elegant short proof we refer to \cite{Jack}.

Given  $(X,\sim)$ as above and $n$ a positive integer we define $(X_n,\sim_n)$  to be the disjoint union of $n$-copies of $X$
\[
X_n=\bigcup_{k=1}^n X\times \{k\}
\]
equipped with the  equivalence relation
\[
(x,j)\sim_n (y,k)\Longleftrightarrow j=k,\;\; x\sim y.
\]
We deduce
\begin{equation}
p(X_n,\sim_n)= \frac{1}{(n|X|)!} \int_0^\infty e^{-x} \bigl(L_\sim(x)\bigr)^n dx
\label{eq: der-n}
\end{equation}
For example, consider the ``marriage relation''
\[
(C,\thicksim),\;\; C=\{\pm 1\},\;\; -1\thicksim 1.
\]
In this case $\hat{C}$ consists of a single  element and the counting function is the number $\nu=2$. Then $(C_n,\thicksim_n)$ can be interpreted as a group of $n$ married couples. If we set
\[
\delta_n:= p(C_n,\thicksim_n)
\]
then we can give the following  amusing interpretation  for $\delta_n$.

\medskip

\noindent {\bf  Couples mixing problem.}  {\sl    At a party  attended  by $n$ couples, the guests were asked to put their names in a hat and then  to  select at random  one name from that pile.  Then the probability  that nobody will select  his/her name or his/her spouse's name is equal to $\delta_n$.  }

\medskip

Using (\ref{eq: der-n}) we deduce
\begin{equation}
\delta_n=\frac{1}{(2n)!}\int_0^\infty e^{-x} \bigl( x^2-4x+2\bigr)^n.
\label{eq: der-c}
\end{equation}
We can ask   about the asymptotic behavior  of the probabilities  $p(X_n,\sim_n)$ as $n\ra \infty$.   In \cite{GIO}, Gillis-Ismail-Offer  describe the first two terms of an asymptotic expansion in powers of $n^{-1}$. To formulate it let us introduce "momenta"
\[
\nu_r=\sum_{c\in \hat{X}}\nu(c)^r.
\]
\begin{theorem}[Gillis-Ismail-Offer]
\begin{equation}
p(X_n,\sim_n)=  \exp( -\frac{\nu_2}{\nu_1})\Bigl(1 -\frac{\nu_1(2\nu_3-\nu_2)-\nu_2^2}{2\nu_1^3}n^{-1}+O(n^{-2}\Bigr)\;\;\mbox{as $n\ra \infty$}.
\label{eq: gio}
\end{equation}
\end{theorem}
For example we deduce from the above that
\begin{equation}
\delta_n =e^{-2} \Bigl( 1-\frac{1}{2} n^{-1} + O(n^{-2}\Bigr),\;\;n\ra \infty
\label{eq: der-ca}
\end{equation}
The proof in \cite{GIO} of  the asymptotic expansion (\ref{eq: gio})  is based on the saddle point technique applied to the integrals in the RHS of (\ref{eq: der-n}) and special properties   of the Laguerre polynomials.

In this paper  we  will investigate the large $n$ asymptotics   of  Laplace transforms
\begin{equation}
\f_n(\q,z)=\frac{z^{dn+1}}{(dn)!}\int_0^\infty e^{-zt}\q(t)^n dt,\;\;\re z >0,
\label{eq: int}
\end{equation}
where $\q(t)$ is a degree $d$ complex  polynomial with leading coefficient $1$.   If we denote  by $\CL[f(t),z]$ the Laplace transform of $f(t)$
\[
\CL[f(t),z]=\int_0^\infty e^{-zt} f(t) dt\;\;\mbox{then}\;\; \f_n(\q,z)=\frac{\CL[\q(t)^n,z]}{\CL[t^{dn},z]}.
\]
The  estimate  (\ref{eq: gio})  will follow from our results by setting
\[
z=1,\;\; \q=  L_\sim.
\]
To formulate the main result we first write $\q$ as a product
\[
\q(t)=\prod_{i=1}^d(t+r_i).
\]
We set
\[
\vec{r} = (r_1,\cdots ,r_d)\in \bC^d,\;\; \mu_s=\mu_s(\vec{r})=\frac{1}{d}\sum_{i=1}^d r_i^s
\]
\begin{theorem}[Existence theorem]    For every $\re z>0$  we have an asymptotic expansion as $n\ra \infty$
\begin{equation}
\f_n(\q,z)= (\sum_{k=0}^\infty  A_k(z) n^{-k}\Bigr).
\label{eq: exp}
\end{equation}
 Above, the  term $A_k(z)$  is  a holomorphic function  on $\bC$  whose coefficients  are universal  elements in the ring of polynomials $\bC(d)[\mu_1,\mu_2,\cdots ,\mu_k]$, where $\bC(d)$ denotes the field    of rational functions in the  variable $d=\deg \q$.
 \label{th: main}
 \end{theorem}

 We can say a bit more about the coefficients $A_k(z)$.

\begin{theorem}[Structure Theorem]  For any $k$ and any degree  $d$ we have
\[
A_k(z)=e^{\mu_1z} B_k(z),
\]
where $B_k\in \bC(d)[\mu_1,\cdots, \mu_k][z]$  is a universal {\bf polynomial} in $z$ with coefficients in  \newline $\bC(d)[\mu_1,\cdots,\mu_k]$.
\label{th: str}
\end{theorem}

 Here is a brief description of the contents. We give the proof of the existence theorem in the next section, while in the  third section we compute the  terms $A_k$ in some cases.  For example  we have
 \[
 A_0(z)=e^{\mu_1z},\;\; A_1(z)= \frac{1}{2d}e^{\mu_1z}(\mu_1^2-\mu_2)z^2,
 \]
 and we can refine (\ref{eq: der-ca}) to
 \begin{equation}
\delta_n =e^{-2} \Bigl( 1-\frac{1}{2} n^{-1} -\frac{23}{96}n^{-2}+ O(n^{-2}\Bigr),\;\;n\ra \infty.
\label{eq: der-ca1}
\end{equation}
These computations  will  lead to a proof of the structure theorem.

 For the reader's convenience    we include  a list of symbols we will use throughout  the paper. The symbol $[n]$ denotes the set $\{1,2,\cdots, n\}$. A $d$-\emph{dimensional (multi)index} will be a  vector $\vec{\alpha}\in \bZ^d_{\geq 0}$.   For  every  vector $\vec{x}\in\bC^d$ and any $d$-dimensional index $\vec{\alpha}$ we   define
\[
\vec{x}^{\vec{\alpha}}=x_1^{\alpha_1}\cdots x_d^{\alpha_d},\;\;|\vec{\alpha}|=\alpha_1+\cdots +\alpha_d,\;\;S(\vec{x})=x_1+\cdots +x_d.
\]
If $n=|\vec{\alpha}|$ then  we define the multinomial coefficient
\[
\binom{n}{\vec{\alpha}}:=\frac{n!}{\prod_{i=1}\alpha_i!}.
\]
These numbers appear in the  \emph{multinomial formula}
\[
S(\vec{x})^n=\sum_{|\vec{\alpha}|=n}\binom{n}{\vec{\alpha}}\vec{x}^{\vec{\alpha}}.
\]

\section{Proof of the existence theorem}
\setcounter{equation}{0}

The key to our approach is the following elementary result.

\begin{lemma} If  $P(x)=p_mt^m+\cdots p_1t +p_0$ is a degree $m$ with complex coefficients then for every  $\re z>1$ we have
\begin{equation}
\frac{\CL[P(t), z]}{\CL[t^m,z]}=\frac{z^{m+1}}{m!}\int_0^\infty e^{-zt} P(t) dt =\sum_{a+b=m} \frac{p_a}{\binom{m}{a}} \frac{z^b}{b!}.
\label{eq: lapl1}
\end{equation}
\label{lemma: lapl}
\end{lemma}

\noindent {\bf Proof}
\[
\frac{z^{m+1}}{m!}\int_0^\infty e^{-zt} P(t) dt= \frac{z^{m+1}}{m!}\sum_{a=0}^mp_a\int_0^\infty e^{-zt}t^a dt=\frac{z^{m+1}}{m!}\sum_{a=0}^mp_a \frac{a!}{z^{a+1}}= \sum_{a+b=m} \frac{p_a}{\binom{m}{a}} \frac{z^b}{b!}.
\]
\qed

Denote by $\q(n,a)$ the coefficient of $t^a$ in $\q(t)^n$.   From (\ref{eq: lapl1}) we deduce
\begin{equation}
\f_n(\q,z)=\sum_{a+b=dn}\frac{\q(n,a)}{\binom{dn}{a}}\frac{z^b}{b!}.
\label{eq: f}
\end{equation}
Using the equality
\[
\q^n=\prod_{i=1}^d\;\underbrace{\Biggl(\sum_{j+k=n}^n\binom{n}{i}t^jr_i^k\Biggr)}_{(t+r_i)^n}
\]
we deduce that if $a+b=dn$ then
\begin{equation}
\q(n,a)=\sum_{|\vec{\alpha}|=b}\;\Biggl(\prod_{i=1}^d\binom{n}{\alpha_j}\Biggr)\, \vec{r}^{\alpha}.
\label{eq: coeff}
\end{equation}
For $|\vec{\alpha}|=b$ we set
\[
B(n,\vec{\alpha}):=\prod_{i=1}^d\binom{n}{\alpha_j},\;\;P_{n,b}(\vec{\alpha}):=\frac{B(n,\vec{\alpha})}{\binom{dn}{|\vec{\alpha}|}},\;\;\rho_b(\vec{\alpha})=\vec{r}^{\vec{\alpha}}.
\]
so that
\begin{equation}
\f_n(\q,z)=\sum_{a+b=dn} \Biggl(\,\sum_{|\vec{\alpha}|=b }P_{n,b} ( \vec{\alpha} )\rho_b(\vec{\alpha} )\,\Biggr)\cdot \frac{z^b}{b!}
\label{eq: f1}
\end{equation}
Observe that  we have
\begin{equation}
P_{n, b}(\vec{\alpha})=\frac{\prod_{i=1}^d\b(1-\frac{1}{n})\cdots (1-\frac{\alpha_i-1}{n})}{\prod_{k=1}^{b-1}(1-\frac{k}{dn})}\cdot\, \underbrace{\frac{1}{d^b} \binom{b}{\vec{\alpha}} }_{:=P_b(\vec{\alpha})}.
\label{eq: distrib}
\end{equation}
The coefficients $P_b(\vec{\alpha})$ define the multinomial probability distribution $P_b$ on the set of multiindices
\[
\Lambda_b=\Bigl\{ \vec{\alpha}\in (\bZ_{\geq 0})^b;\;\;|\vec{\alpha}|=b\Bigr\}
\]
For every random variable $\zeta$ on $\Lambda_b$ we denote by $E_b(\zeta)$ its expectation with respect to the probability distribution $P_b$.  For each $n$ we have  a random variable $\zeta_{n,b}$ on $\Lambda_b$ defined by
\[
\zeta_{n,b}(\vec{\alpha})= \frac{\prod_{i=1}^d\b(1-\frac{1}{n})\cdots (1-\frac{\alpha_i-1}{n})}{\prod_{k=1}^{b-1}(1-\frac{k}{dn})}\rho_b(\vec{\alpha}).
\]
Form (\ref{eq: f1}) and (\ref{eq: distrib}) we deduce
\begin{equation}
\f_n(\q,z)=\sum_{a+b=dn}E_b(\zeta_{n,b}) \frac{z^b}{b!}.
\label{eq: f2}
\end{equation}
To find the asymptotic expansion for $\f_n$ we will find  asymptotic expansions  in powers of $n^{-1}$  for the expectations $E_b(\zeta_{n,b})$ and them add them up using (\ref{eq: f2}).

 For every nonnegative integer  $\alpha$ we   define a polynomial
\[
W_\alpha(x)=\left\{
\begin{array}{ccc}
1 &{\rm if} & \alpha =0,1\\
&&\\
\prod_{j=1}^{\alpha-1}(1-jx) & {\rm if} & \alpha >1
\end{array}
\right.
\]
For a $d$-dimensional multi-index     $\vec{\alpha}$ we set
\[
W_{\vec{\alpha}}(x)=\prod_{i=1}^d W_{\alpha_i}(x).
\]
We can now rewrite  (\ref{eq: distrib}) as
\[
P_{n,b}(\vec{\alpha}) = P_b(\vec{\alpha}) \frac{W_{\vec{\alpha}}(\frac{1}{n})}{W_b(\frac{1}{dn})}.
\]
We set
\[
R_b(\vec{\alpha},x)={W_{\vec{\alpha}}(x)},\;\;K_b(\vec{\alpha},x) =\frac{1}{W_b(\frac{x}{d})} R_b(\vec{\alpha},x)\rho_b(\alpha).
\]
We regard  the correspondences
\[
\vec{\alpha}\mapsto R_b(\vec{\alpha},x),\; K_b(\vec{\alpha},x)
\]
as a  random variables  $R_b(x)$ and $K_b(x)$ on $\Lambda_b$ valued in the field of  rational functions. We deduce
\[
\zeta_{n,b}=K_b(n^{-1}).
\]
Observe
\[
E_b(x) =E_b(K_b(x))=\frac{1}{W_b(x)} E_b(R_b(x))
\]
From the fundamental  theorem of symmetric polynomials we deduce  that the expectations   $E_b(R_b(x))$ are \emph{universal} polynomials
\[
E_b(R_b(x))\in \bC[\mu_1,\cdots,\mu_b][x], \;\; \deg_x E_b(R_b(x))\leq b-d,
\]
whose coefficients  have  degree $b$ in the variables $\mu_i$, $\deg \mu_i=i$.  We deduce that  $E_b(x)$ has a Taylor series  expansion
\[
E_b(x) =\sum_{m\geq 0} E_b(m) x^m
\]
such that $E_b(m)\in \bC(d)[\mu_1,\cdots, \mu_b]$. The rational function $x\ra K_b(\vec{\alpha},x)$ has a Taylor expansion   at $x=0$ convergent for $|x|<\frac{d}{b-1}$ so the above series  converges for $|x|<\frac{d}{b-1}$. We would like to estimate the size of the coefficients $E_b(m)$.  The tricky part is that the radius of convergence of $E_b(x)$ goes to zero as $b\ra \infty$.

\begin{lemma}   Set
\[
R=\max_{1\leq i\leq d}|r_i|.
\]
There exists a constant $C$ which depends only on $R$ and $d$ such that for every  $b\geq 0$ and every  $1\leq \lambda_b <\frac{b}{b-1}$  we have the inequality
\begin{equation}
|E_b(m)|= \Bigl(\frac{b}{\lambda_b d}\Bigr)^{m} C^b\frac{b^{b-1}}{(b-2)!(1-\lambda_b\frac{b-1}{b})}.
\label{eq: coef-est}
\end{equation}
\label{lemma: coef-est}
\end{lemma}

\noindent {\bf Proof}\hspace{.3cm}   Note first that
\[
|\rho_b(\vec{\alpha})|\leq R^b,\;\;\forall|\vec{\alpha}|=b.
\]
For $b=0,1$ we deduce form the definition of the polynomials $W_\alpha$ that $E_b(x)=1$.  Fix $m$ and $b>1$. Using the Cauchy residue formula we deduce
\[
E_b(m)=\frac{1}{2\pi\sqrt{-1}}\int_{|x|= \hbar}\frac{1}{x^{m+1}} E_b(x) dx,\;\;\hbar= \lambda_b\,\cdot\, \frac{d}{b}.
\]
Hence
\[
|E_b(m)|\leq \frac{1}{\hbar^{m}}\sup_{|x|=\hbar} |E_b(x)|\leq \frac{b^{m}R^b}{(\lambda_bd)^{m} \min_{|x|=\hbar}|W_b(x/d)|}\,\cdot\, \max_{|x|=\hbar} E_b(R_b(x)).
\]
Next observe that
\[
W_b(x/d)= (b-1)!\prod_{k=1}^{b-1}(\frac{1}{k} -x/d)),\;\; \hbar/d<1/k,\forall k \leq b-1,
\]
from which we conclude
\[
\min_{|x|=\hbar}|W_b(x)|= W_b(\hbar)=\prod_{k=1}^{b-1}(1 -\frac{k\lambda_b}{b})=\frac{1}{b^{b-1}}\prod_{k=1}^{b-1}(b-k\lambda_b)\geq \frac{(b-2)!(1-\lambda_b\frac{b-1}{b})}{b^{b-1}}.
\]
To estimate $E_b(R_b(x))$ from above observe that for every $1\leq k\leq (b-1)$  and $|x|=\hbar$ we have
\[
|1-kx|\leq  1+ k|x|=1+\frac{k\lambda_bd}{b}< 1+ d
\]
This shows that for every $|\vec{\alpha}|=b$ and $|x|=\hbar$ we have
\[
|R_b(\vec{\alpha},x)|< (1+d)^b.
\]
The lemma follows by assembling all the facts established above.

\qed

Define the formal power series
\[
A_m(z):=\sum_{b\geq 0}  E_b(m) \frac{z^b}{b!}\in \bC\{z\}.
\]
 The estimate  (\ref{eq: coef-est}) shows that this series converges for all $z$.

For every   formal power series  $f=\sum_{k\geq 0} a_kT^k$   and every  nonnegative integer $\ell$ we denote by  $J_T^\ell(f)$ its $\ell$-th jet
\[
J_T^\ell(f)=\sum_{k=0}^\ell a_k T^k.
\]
 For $x=n^{-1}$ we have
\[
\f_x(z)=\f_n(\q,z)=\sum_{b\leq d/x}  E_b(x)\frac{z^b}{b!}=\sum_{b\leq d/x}  \Bigl(\sum_{m\geq 0}E_b(m)x^m\Bigr)\frac{z^b}{b!}
\]
\[
=\sum_{m\geq 0}\Bigl(\sum_{b\leq d/x}E_b(m)\frac{z^b}{b!}\Bigl)x^m=\sum_{m\geq 0} J^{d/x}_z (A_m(z))x^m.
\]
Consider the formal power series in $x$ with coefficients in the ring $\bC\{z\}$ of convergent power series in $z$
\[
\f_\infty(z)=\sum_{m\geq 0} A_m(z) x^m\in \bC\{z\}[[x]].
\]
We will prove that for every $\ell \geq 0$  and every $z\in \bC$ we have
\begin{equation}
|\f_n(z)- J^\ell_x \f_\infty(z)| =O(n^{-\ell-1}),\;\;\mbox{as $n\ra \infty$}.
\label{eq: ultimate}
\end{equation}
To prove this  it is convenient to introduce the ``rectangles''
\[
D_{u,v}=\Bigl\{ (b ,m)\in (\bZ_{\geq 0})^2;\;\; b\leq u,\;\;  m\leq v.\Bigr\}
\]
In this notation we have ($x=n^{-1}$)
\[
\f_n(z)=\sum_{(b,m)\in D_{n,\infty}}E_b(m) x^m\frac{z^b}{b!},\;\; J^\ell_x \f_\infty(z)=\sum_{(b,m)\in D_{\infty,\ell}}E_b(m) x^m\frac{z^b}{b!}.
\]
Then
\[
\f_n(z)- J^\ell_x \f_\infty(z) =\underbrace{\sum_{b\leq dn}\Bigl(\sum_{m> \ell} E_b(m) x^m\Bigr)\frac{z^b}{b!}}_{S_1(n)} \, +\, \underbrace{\sum_{m\leq \ell}\Bigl(\sum_{b> dn} E_b(m) \frac{z^b}{b!}\Bigr)x^m}_{S_2(n)}.
\]
We estimate each    sum separately. Using  (\ref{eq: coef-est}) with  a $\lambda_b>1$ to be  specified later we deduce
\[
\Bigl|\sum_{m> \ell} |E_b(m) x^m|\leq \frac{C^bb^{b-1}}{(b-2)!(1-\lambda_b\frac{b-1}{b})}\sum_{m >\ell} \Bigl(\frac{bx}{\lambda_bd}\Bigr)^m.
\]
The inequality $b\leq dn$ can be translated into $\frac{bx}{d}\leq 1$ so that the above series is convergent  for $b\leq dn$   whenever $\lambda_b> 1$ so that
\[
\Bigl|\sum_{m> \ell} |E_b(m) x^m|\leq  \frac{C^bb^{b-1}}{(b-2)!(1-\lambda_b\frac{b-1}{b})}\Bigl(\frac{bx}{\lambda_bd}\Bigr)^{\ell+1}\frac{1}{1-\frac{bx}{\lambda_b d}}.
\]
When $b\leq dn$ we have
\[
1-\frac{bx}{\lambda_b d} >1-\frac{1}{\lambda_b}.
\]
If we choose
\[
\lambda_b=\Bigl( \frac{b}{b-1}\Bigr)^{1/2}
\]
we deduce
\[
1-\lambda_b\frac{b-1}{b})=1-\Bigl(\frac{b-1}{b}\Bigr)^{1/2} \Longrightarrow\frac{1}{1-\lambda_b\frac{b-1}{b}} <b\
\]
and, since $\frac{bx}{\lambda_b d}\leq \frac{b}{d}x$,
\[
\frac{1}{1-\frac{bx}{\lambda_b d}}<\frac{1}{1-\frac{1}{\lambda_b}} <2b.
\]
Using the inequalities
\[
k! >\Bigl(\frac{k}{e}\Bigr)^k,\;\;\forall k>0
\]
we conclude that  for $b\leq dn$ we have
\[
\sum_{m> \ell} |E_b(m) x^m|\leq C_1^bb^{\ell+2}x^{\ell+1}  .
\]
Since the series $\sum_{b\geq 0} {C_1^bb^{\ell+2}}\frac{z^b}{b!}$ converges we conclude that
\[
S_1(n)= O(x^{\ell+1}).
\]
To estimate the second sum we choose $\lambda_b=1$ in  (\ref{eq: coef-est}) and we deduce
\[
E_b(m)\leq  C_3^b.
\]
Hence
\[
\Bigl|\sum_{b> dn} E_b(m) \frac{z^b}{b!}|\leq \frac{(C_3|z|)^b b^2}{b!} < (2C_3|z|)^2 \sum_{b>dn}\frac{(|C_3|z|)^{b-2}}{(b-2)!}
\]
Using  Stirling's formula we deduce that for fixed $z$ we have
\[
\sum_{b>dn} \frac{(|C_3|z|)^{b-2}}{(b-2)!} < C_4(z) n^{-\ell-1}.
\]
Hence
\[
|S_2(n)| \leq C_4(z)(\ell+1) n^{-\ell-1}.
\]
This completes the proof of  (\ref{eq: ultimate}) and of  Theorem \ref{th: main}.

\qed

\section{Additional structural results}
\setcounter{equation}{0}

\subsection{The case $d=1$.}  Hence
\[
\q(t)=(t+\mu_1).
\]
We have
\[
\int_0^\infty e^{-zt} (t+\mu_1)^n dt=  e^{\mu_1 z}\int_0^\infty e^{-zt}t^n dt = e^{\mu_1z} \frac{n!}{z^{n+1}}
\]
Hence in this case
\[
\f_n(z)= e^{\mu_1 z}
\]
and we deduce
\[
A_0(z)=e^{\mu_1z},\;\;A_k(z)=0 ,\;\;\forall k\geq 1.
\]

\subsection{The case  $d=2$.} This  is a bit  more complicated.     We assume first that $\mu_1=0$ so that
\[
\q(t)= t^2-\si^2
\]
Then
\[
\q(n,a)= \left\{
\begin{array}{ccc}
(-1)^k\si^{2(n-k)} \binom{n}{k} & {\rm if } & a=2k\\
& &\\
0   & {\rm if} &   \mbox{$a$ is odd}
\end{array}
\right. ,
\]
and we deduce
\[
\f_n(z)=\sum_{b=0}^n\frac{(-1)^{b} \binom{n}{n-b} }{\binom{2n}{2n-2b}} \frac{(\si z)^{2b}}{(2b)!}= \sum_{b=0}^n \frac{n!(2n-2b)!}{(n-b)!(2n)!} \frac{(-1)^b(\si z)^{2b}}{b!}
\]
\[
=\sum_{b=0}^n\frac{n(n-1)\cdots (n-b+1)}{2n(2n-1)\cdots (2n-2b+1)}\frac{(-1)^b(\si z)^{2b}}{b!}
\]
\[
=\sum_{b=0}^n\frac{1}{2^{2b}}n^{-b}\frac{(1-1/n)\cdots (1-(b-1)/n)}{(1-1/(2n)\cdots (1- (2b-1)/(2n)}\frac{(-1)^b(\si z)^{2b}}{b!}.
\]
\[
=1-\frac{1}{2}n^{-1}\frac{1}{1-1/(2n)}\frac{(\si z^)2}{2!}+ \frac{1}{2^4}n^{-2}\frac{(1-1/n)}{(1-1/(2n))(1-2/(2n))(1-3/(2n))}\frac{(\si z)^4}{4!} +\cdots
\]
To obtain $A_k(z)$ we need to collect the powers $n^{-k}$.   The above description shows that the coefficients of the monomials $z^{2b}$  contain only powers $n^{-k}$, $k\geq b$. We conclude that  $A_k(z)$ is a polynomial and
\[
\deg_z A_k(z)\leq 2k.
\]
Let us compute the first few of these polynomials. We have
\[
\f_n(z)=1-\frac{1}{2}n^{-1}\Bigl(1+\frac{1}{2}n^{-1}+\cdots\Bigr)\frac{(\si z)^2}{2!}+ \frac{1}{2^4}n^{-2}\Bigl(1+\cdots\Bigr)\frac{(\si z)^4}{4!}+\cdots
\]
We deduce
\[
A_0(z)=1,\;\;A_1(z)=-\frac{1}{4}(\si z)^2,\;\;A_2(z)= -\frac{1}{8}(\si z)^2+\frac{1}{2^44!}(\si z)^4.
\]
If $\mu_1\neq 0$  so that
\[
\q(t)=(t+r_1)(t+r_2),\;\;r_1+r_2=2\mu_1,
\]
then we make the change in variables  $t=s-\mu_1$ so that
\[
\q(t)= P(s)=s^2-r^2,\;\;\si^2=(r_1-\mu_1)^2=    \frac{(r_1-r_2)^2}{4}.
\]
Now observe that
\[
4\mu_1^2+(r_1-r_2)^2=(r_1+r_2)^2+(r_1-r_2)^2= 2(r_1^2+r_2^2)= 4\mu_2
\]
so that
\[
\si^2= \mu_2-\mu_1^2.
\]
Then
\[
\f_n(\q,z)=\frac{z^{2n+1}}{(2n)!}\int_0^\infty e^{-zt} \q(t)^n= \frac{z^{2n+1}}{(2n)!}\int_0^\infty e^{-z(s-\mu_1)} P(s)^n ds=e^{\mu_1z} \f_n(P,z).
\]
We deduce
\begin{equation}
A_0(z)=e^{\mu_1 z},\;\; A_1(z) =-\frac{e^{\mu_1 z}}{4}(\si z)^2,\;\; A_2(z)=e^{\mu_1 z}\Bigl(-\frac{1}{8}(\si z)^2+\frac{1}{2^44!}(\si z)^4\Bigr).
\label{eq: asy1}
\end{equation}
For the  couples mixing problem we have
\[
\q(t)= t^2-4t+2
\]
so that
\[
\mu_1=-\frac{4}{2}=-2,\;\; \si^2=\frac{1}{4}(r_1-r_2)^2= \frac{1}{4}\Bigl( (r_1+r_2)^2-4r_1r_2\Bigr)=\frac{1}{4}( 16-8)= 2.
\]
and we deduce
\begin{equation}
\delta_n= \f_n(\q,z=1)=e^{-2}\Bigl(1-\frac{1}{2}n^{-1}-\frac{23}{96}n^{-2}+O(n^{-3})\Bigr).
\label{eq: mix}
\end{equation}

\subsection{The general case.} Let us determine the coefficients $A_0(z)$ and $A_1(z)$ for general degree $d$. We use the definition
\[
A_k(z)=\sum_{b\geq 0} E_b(k)\frac{z^b}{b!}.
\]
For $|\vec{\alpha}|=b$
\[
W_{\vec{\alpha}}(x)=W_{b,\alpha}(x)=\prod_{i=1}^d \Bigl(\prod_{j=1}^{\alpha_i-1} (1-jx)\Bigl)
\]
\[
=\prod_{i=1}^d\Biggl(1 -\Bigl(\sum_{j=1}^{\alpha_i-1}j\Bigr) x+\cdots\Biggr)=1-\frac{1}{2}\Bigl(\sum_{i=1}^d\alpha_i(\alpha_i-1)\Bigr)x+\cdots .
\]
\[
W_b(x/d)= \prod_{k=1}^{b-1}(1+jx/d+\cdots)= 1+\frac{b(b-1)}{2d} x+\cdots .
\]
Next   compute the expectation of $R_b(x)$
\[
E_b(R_b(x))= E_b(\rho_b) - \frac{1}{2}E_b\Bigl(\sum_{i=1}^d\alpha_i(\alpha_i-1)\vec{r}^{\vec{\alpha}}\Bigr)x+\cdots.
\]
The multinomial formula implies
\[
E_b(\rho_b)= \mu_1^b.
\]
Next
\[
E_b\Bigl(\sum_{i=1}^d\alpha_i(\alpha_i-1)\vec{r}^{\vec{\alpha}}\Bigr)=\frac{1}{d^b}\sum_{|\vec{\alpha}|=b}\binom{b}{\vec{\alpha}}(\sum_{i=1}^d\alpha_i(\alpha_i-1))\vec{r}^{\vec{\alpha}}
\]
Now consider the partial differential operator
\[
\p=\sum_{i=1}^d r_i^2 \frac{\pa^2}{\pa r_i^2}.
\]
Observe that the monomials $\vec{r}^{\vec{\alpha}}$ are eigenvectors of $\p$
\[
\p\vec{r}^{\vec{\alpha}}= (\sum_{i=1}^d\alpha_i(\alpha_i-1))\vec{r}^{\vec{\alpha}}.
\]
We deduce
\[
E_b\Bigl(\sum_{i=1}^d\alpha_i(\alpha_i-1)\vec{r}^{\vec{\alpha}}\Bigr)=\frac{1}{2d^b}\p S(\vec{r})^b=\frac{1}{2}\p\mu_1^b.
\]
Hence
\[
E_b(R_b(x)=\mu_1^b-\frac{1}{2}(\p\mu_1^b)x+\cdots
\]
and we deduce
\[
E_b(x)=(\mu_1^b-\frac{1}{2}(\p\mu_1^b)x+\cdots)(1+\frac{b(b-1)}{2d} x+\cdots )
\]
\[
=\mu_1^b+\frac{1}{2}\Bigl(\frac{b(b-1)}{d}\mu_1^b-\p\mu_1^b\Bigr) x+\cdots
\]
We deduce $A_0(z)= e^{\mu_1 z}$
\[
 A_1(z)=\frac{\mu_1^2}{2d}\sum_{b=2}^\infty \frac{z^b}{(b-2)!}-\frac{1}{2}\p e^{\mu_1z}=\frac{\mu_1^2z^2}{2d} e^{\mu_1z}-\frac{1}{2}\p e^{\mu_1 z}.
\]
We can simplify the answer some more.
\[
\p\mu_1^b=\frac{1}{d^b}\p S(x)^b=\frac{b(b-1)}{d^b} \Bigl(\sum_{i=1}^d r_i^2\Bigr) S(x)^{b-2}=\frac{b(b-1)}{d}\mu_2\mu_1^{b-2}.
\]
We conclude that
\[
\p e^{\mu_1 z}=\frac{\mu_2z^2}{d}\sum_{b\geq 2}\frac{ (\mu_1z)^{b-2} }{(b-2)!}= \frac{\mu_2z^2}{d} e^{\mu_1z}.
\]
Hence
\begin{equation}
A_0(z)=e^{\mu_1z},\;\;A_1(z)=\frac{e^{\mu_1z}}{2d}(\mu_1^2-\mu_2)z^2
\label{eq: asy2}
\end{equation}
For $d=2$ we recover part of the formul{\ae} (\ref{eq: asy1}).

\medskip

\subsection{The proof of Theorem \ref{th: str}} Clearly we can assume $d>1$. We imitate the strategy  used in the case $d=2$. Thus, after the change in variables $t\ra t-\mu_1$ we can assume that $\mu_1=0$ so that  $\q(t)$ has the special form
\[
\q(t)=t^d+ a_{d-2} t^{d-2} +\cdots + a_0.
\]
Set
\[
T(n,b):=\frac{\q(n, dn-b)}{\binom{dn}{dn-b}}
\]
This is a power series in $x=n^{-1}$,
\[
T(n,b)= T_b(x)\!\mid_{x=n^{-1}},\;\; T_b(x)=\sum_{k\geq 0} T_b(k) x^k.
\]
We have
\[
A_k(z)=\sum_{b\geq 0} T_b(k)\frac{z^b}{b!},
\]
 and we need to prove that $A_k$ is a polynomial for every $k$. We denote by $\ell(b)$ the order of of the first non-zero  coefficient of $T_b(x)$,
 \[
 \ell(b)=\min\{k\geq 0;\;\; T_b(k)\neq 0\}.
 \]
To prove the desired conclusion it suffices to show that
\begin{equation}
\lim_{b\ra \infty} \ell(b)=\infty.
\label{eq: lim}
\end{equation}
For every multiindex $\vec{\beta}=(\beta_d,\beta_{d-2},\cdots,\beta_1,\beta_0)$ we set
\[
L(\vec{\beta})= d\beta_d +(d-2)\beta_{d-2}+\cdots +\beta_1.
\]
We set $\vec{a}= (1,a_{d-2},\cdots , a_1, a_0)\in {\bC}^d$. We have
\begin{equation}
T(n,b)=\frac{1}{\binom{dn}{dn-b}}\, \cdot \, \sum_{|\vec{\beta}|=n, L(\vec{\beta})=dn-b}\binom{n}{\vec{\beta}} \vec{a}^{\vec{\beta}}.
\label{eq: sum}
\end{equation}
Now observe that  we have
\[
 d|\vec{\beta}-L(\vec{\beta})={2}\beta_{d-2}+{3}\beta_{d-3}+\cdots (d-1)\beta_1 + d\beta_0 =d|\vec{\beta}|-L(\vec{\beta})= b.
 \]
In particular we deduce
\begin{equation}
\beta_j\leq \frac{b}{d-j}\leq \frac{b}{2},\;\;\forall 0\leq j\leq d-2
\label{eq: cruc1}
\end{equation}
and
\[
2\beta_d +b =2\beta_d={2}\beta_{d-2}+{3}\beta_{d-3}+\cdots (d-1)\beta_1 + d\beta_0
\]
\[
\geq  2\beta_d+2\beta_{d-2} +\cdots + 2\beta_1+2\beta_0 = 2n
\]
so that
\begin{equation}
n-\beta_d\leq \frac{b}{2}.
\label{eq: cruc2}
\end{equation}
These simple observations have several important consequences.

First, observe that they imply  that there exists  an integer $N(b)$ \emph{which depends only $b$ and $d$}, such that for  any $n>0$ the  number of multi-indices $\vec{\beta}$ satisfying  $|\vec{\beta}|=n$ and $L(\vec{\beta} )= dn-b$ is $\leq N(b)$. Thus the sum (\ref{eq: sum}) has fewer than $N(b)$ terms.

Next, if  we set $|a|:=\max_{0\leq_j\leq d-2}|a_j|$ then, we deduce
\[
|\vec{a}^{\vec{\beta}}|\leq |a|^{\beta_0+\cdots +\beta_{d-2}} \leq|a|^{\frac{b(d-1)}{2}} =  C_5(b).
\]
Finally, using the identity
\[
\binom{n}{\vec{\beta}}= \binom{n}{\beta_d}\cdot \binom{n-\beta_d}{\beta_{d-2}}\binom{n-\beta_d-\beta_{d-2}}{\beta_{d-3}}\cdots
\]
 the inequalities (\ref{eq: cruc2})  and  $\binom{m}{k}\leq 2^m$, $\forall m\geq k$ we deduce
 \[
\binom{n}{\vec{\beta}}\leq \binom{n}{\beta_d}\cdot 2^{\frac{b(d-1)}{2}}\leq 2^{\frac{b(d-1)}{2}}\binom{n}{\lfloor b/2\rfloor +1} \leq C_6(b) n^{\lfloor b/2\rfloor +1},\;\;\forall n\gg b.
\]
Hence
\[
\sum_{|\vec{\beta}|=n, L(\vec{\beta})=dn-b}\Bigl|\binom{n}{\vec{\beta}} \vec{a}^{\vec{\beta}}\Bigr|\leq N(b) C_5(b)C_6(b) n^{\lfloor b/2\rfloor +1} =C_7(b) n^{\lfloor b/2\rfloor +1}.
\]
On the other hand
\[
\frac{1}{\binom{dn}{dn-b}}\leq C_8(b) n^{-b}
\]
so that
\[
|T(n,b)|= |T_b(n^{-1})|\leq C_9(b) n^{\lfloor b/2\rfloor+1-b} \leq C_9(b) n^{1-b/2}.
\]
This shows
\[
T_b(k)=0,\;\;\forall k\leq b/2-1
\]
so that
\[
\ell(b)\geq b/2-1\ra \infty\;\;\mbox{as $b\ra \infty$}.
\]
\qed

\begin{remark} We can say a bit more  about the structure of  the polynomials
\[
B_k(\mu_1,\cdots, \mu_d,z)\in R_d=\bC[\mu_1,\cdots,\mu_d,z], \;\; k>0.
\]
If we regard $B$ as a polynomial in $r_1,\cdots,r_d$ we see that it vanishes precisely when $r_1=\cdots=r_d$.   Note that
\[
r_1=\cdots =r_d=r\Longleftrightarrow \q(t)= (t+r)^d.
\]
On the other hand
\[
\sum_k t^k\mu_k= \frac{1}{d}\sum_{i=1}^d\sum_{k\geq 0}(r_it)^k=  \frac{1}{d}\sum_{i=1}^d\frac{1}{1-r_it} \stackrel{(s:= 1/t)}= \frac{s}{d}\sum_{i=1}^d\frac{1}{s+\mu_i}=\frac{s}{d} \frac{\q'(s)}{\q(s)}.
\]
If $\q(s)=(s+r)^d$ we deduce
\[
\frac{s}{d} \frac{\q'(s)}{\q(s)}=\frac{s}{s+r}=\frac{1}{1-rt}=\sum_{k\geq 0} (rt)^k.
\]
This implies that
\[
r_1=\cdots =r_d \Longleftrightarrow \mu_i^j=\mu_j^i,\;\;\forall 1\leq i,j\leq k\Longleftrightarrow \mu_j=\mu_1^j,\;\;\forall 1\leq j\leq d.
\]
The ideal $I$ in $R_d$ generated by the binomials $\mu_1^j-\mu_j$ is prime since  $R_d/I\cong \bC[\mu_1,z]$.  Using the Hilbert \emph{Nullstelens{a}tz} we  deduce  that $B_k$ must belong to this ideal so that we can write
\[
B_k(\mu_1,\cdots,\mu_d,z)=  A_{2k}(\mu,z)(\mu_1^2-\mu_2)+\cdots +A_{dk}(\mu,z)(\mu_1^d-\mu_d).
\]
\qed
\end{remark}

\end{document}